\documentclass[10pt,notitlepage,twoside,a4paper]{amsart}
 \usepackage{amsfonts}

\usepackage{amsmath,amssymb,enumerate}

\usepackage{epsfig,fancyhdr,color}

\usepackage{amssymb}
\usepackage{amsmath,amsthm}
\usepackage{latexsym}
\usepackage{amscd}
\usepackage{psfrag}
\usepackage{graphicx}
\usepackage[latin1]{inputenc}
\usepackage[all]{xy}
\usepackage[mathcal]{eucal}

\definecolor{NoteColor}{rgb}{1,0,0}

\renewcommand{\textsc}{\textcolor{red}}

\newtheorem*{theorem 1}{\rm\bf Proposition 1}
\newtheorem*{theorem 2}{\rm\bf Proposition 2}

\theoremstyle{definition}

\theoremstyle{remark}

\def\interieur#1{\mathord{\mathop{\kern 0pt #1}\limits^\circ}}

\title[Euler and Chebyshev]{Euler and Chebyshev:
\\ From the sphere to the plane and backwards}

\author[Athanase Papadopoulos]{Athanase Papadopoulos}

 \address{Athanase Papadopoulos, Institut de Recherche Math\'ematique Avanc\'ee\\ CNRS et Universit\'e de Strasbourg\\\small 7 rue Ren\'e
  Descartes - 67084 Strasbourg Cedex, France}

\date{\today}

\begin{document}

  \maketitle
    \begin{abstract}
 We report on the works of Euler and Chebyshev on the drawing of geographical maps. We point out relations with questions about the fitting of garments that were studied by Chebyshev.

   \end{abstract}
\medskip
\medskip
\medskip

 This paper will appear in the Proceedings in Cybernetics, 
 a volume dedicated to the 70th anniversary of Academician Vladimir Betelin.
     
 \medskip
      \medskip

 Keywords: Chebyshev, Euler, surfaces, conformal mappings, cartography, fitting of garments, linkages.

 AMS classification: 30C20,  91D20, 01A55, 01A50, 53-03, 53-02, 53A05,  53C42, 53A25.

  \section{Introduction}
  
  Euler and Chebyshev were both interested in almost all problems in pure and applied mathematics and in engineering, including the conception of industrial machines and technological devices.
In this paper, we report on the problem of drawing geographical maps on which they both worked. It is not surprising that the two men were attracted by this problem which involves nontrivial questions of geometry and analysis, and which at the same time has practical applications. Euler and Chebyshev brought into this field all their knowledge from differential geometry, the theory of differential equations, and the calculus of variations. At the same time, cartography raised new theoretical questions, in geometry and analysis. One may also note here that Chebyshev was particularly interested in approximation and interpolation theories, and in particular interpolation using the method of least squares, and he used these theories in his work on cartography. In what follows, we shall describe and compare the works of Euler and Chebyshev on cartography, making relations with other problems on which they worked.

From the mathematical point of view, the problem of drawing geographical maps is the one of mapping on a Euclidean plane a subset of a curved surface, which is usually a sphere, or a spheroid, representing the Earth\footnote{It was known since the work of Newton that the Earth is spheroidal and not spherical, namely, it is slightly flattened at the poles. 
    Of course, at the local level, the Earth can be neither spherical nor spheroidal; it has mountains, valleys, canyons, etc. All this has to be taken into account in the drawing of maps of very small regions.} or the Celestial sphere. It was already known to the mathematicians and geographers of Greek antiquity (although at that time there was no mathematical proof of this fact) that a map from a region of a sphere onto the Euclidean plane cannot preserve at the same time angles, ratios of distances and ratios of areas.  The question then was to find maps which preserve ``in the best way" angles,  distances or areas, or an appropriate compromise between these quantities. This theme led to many interesting mathematical developments. Several prominent mathematicians worked on it from various points of view.  One may mention here Ptolemy, the famous  mathematician, astronomer and geographer from the second century A.D. who thoroughly studied these questions. From the modern period, besides Euler and Chebyshev, one should mention at least the names of  Lambert, Lagrange, Gauss, Beltrami and Darboux, and there are many others.

The second problem which we consider in this report was studied by Chebyshev. This is the problem of binding a surface with a piece of fabric. The piece of fabric in which we are interested is made of two families of threads that are perpendicular to each other, forming a Euclidean net of rectangles. A garment is made with this piece of fabric, and this garment is intended to cover a surface, which is usually a part of the human body, and optimally it should take its form.  For that, the small rectangles formed by the net are deformed. During this deformation the lengths of the sides of the rectangles remain constant, but the angles they make change. The cloth becomes a curved surface, and the rectangles become curvilinear parallelograms, that is, quadrilaterals in space bounded by two pairs of equal sides.

  The problem of drawing geographical maps and the one of fitting of garments are inverse of each other. Indeed, on the one hand, one searches for mappings from a piece of the sphere into a Euclidean piece of paper with minimal distortion, and on the other hand, one constructs a map from a Euclidean piece of fabric onto a curved surface (part of a human body), such that the fabric fits the curved surface  with minimal distortion.
  
  In this paper, we consider these two problems in the way they were addressed by Euler and Chebyshev. We point out relations between these problems and other works of  Euler, Chebyshev and some other scientists.
  
  \section{Euler and Chebyshev}
  
  We start with a few words on Euler.

   Leonhard Euler published memoirs and books on all the fields of pure and applied mathematics that were known at his epoch: geometry (Euclidean and spherical), analysis, differential equations, number theory, probability, cartography, astronomy, etc. He created the field of calculus of variations. He is also one of the founders of topology.
   
  Euler is the father of Russian mathematics. During his two long stays at the Academy of Sciences of Saint Petersburg (1727---1741 and 1766---1783) he contributed to the training of Russian students in mathematics.\footnote{Euler's Russian students include  Semion Kirillovich Kotelnikov, Stepan Rumovsky, Ivan Iudin, Petr  Inokhodtsov, Mikhail Evseyevich Golovine, and there are many others. It is known that Euler liked his Russian students. During his stay in Berlin, several of them followed him. They also helped him in translating his books into Russian. Romovsky translated Euler's \emph{Letters to a German Princess}. Iudin and Inokhodtsov translated Euler's \emph{Vollst\"andige Anleitung zur Algebra} (Complete introduction to algebra). Golovine became professor at Saint Petersburg's pedagogical school and he translated into Russian Euler's \emph{Th\'eorie compl\`ete de la construction et de la man\oe uvre des vaisseaux} \cite{E426}.} He wrote for the young generations textbooks that were used in universities, in military schools, in naval institutes and other technical schools. He also wrote on basic school arithmetic, \emph{Einleitung zur Rechen-Kunst zum Gebrauch des Gymnasii bey der Kayserlichen Academie der Wissenschafften in St. Petersburg.} (Introduction to the Art of Reckoning, for use in the Gymnasiums of the Imperial Academy of Sciences in Saint Petersburg) \cite{E17}.
   
    Euler was also a physicist. He published on mechanics, acoustics, optics, fluid dynamics, astronomy, geodesy and other fields of physics.  From the practical point of view, Euler conducted astronomical and acoustical measures, and he conceived new optical instruments. In Saint Petersburg, he solved problems concerning floods. 
 During his stay in Berlin (1741--1766), he was solicited by his mentor, the king Frederick II, to solve various concrete problems. This is testified in the correspondence between the two men.\footnote{Euler's Opera Omnia, Ser IV1, vol. VI, is devoted to the correspondence between Euler, Maupertuis and Frederick II.} For instance, in a letter dated April 30, 1749, the monarch asks Euler's help on problems concerning the navigation between the two rivers Oder and Havel. During the Seven Years' War (1756--1763), he wrote to Euler to ask him to improve the technique of the spyglass. Euler's \emph{Opera Omnia} contain at least 
 22 papers on the theory of machines. We mention as examples a few of these papers, whose titles are significant, \emph{De machinis in genere}
(On machines in general) \cite{E194} (1753), \emph{Sur l'action des scies} 
(On the action of saws) \cite{E235} (1758), \emph{Sur la force des colonnes}
(On the strength of columns) \cite{E238} and \emph{Recherches plus exactes sur l'effet des moulins \`a vent} 
(More exact researches on the effect of windmills) \cite{E233} (1758). There are also several memoirs by Euler on the technique and conception of optical instruments, e.g. 
 \cite{E118}, \cite{E196} and \cite{E239}. Euler also published a translation with commentary of a famous book on gunnery by 
 Benjamin Robins \cite{E77}, with the German title 
 \emph{Neue Grunds\"atze der Artillerie} (New principles of artillery). In this book Euler investigates, among other questions, the nature of air and fire and the motion of bodies that are projected in the air. He  establishes relations between elasticity, density and temperature.  
 Euler also published several books and papers on ship building, e.g. \cite{Euler-Navalis}, \cite{E426} and \cite{E520}. This subject was very important to the Russian rulers, who were at the same time the patrons of the Academy of Sciences of Saint Petersburg. The extension of the Russian fleet was one of their dearest projects. As a matter of fact, this was one of the major concerns of Peter the Great, when he planned the Russian Academy of Sciences. He was aware of the role of mathematics and mathematicians in the training of teachers in naval, military and engineering schools. In any case, Euler was most interested in naval construction, a rich subject which involves hydrostatics, architecture, resistance, motion, machines, and stability theory. Euler also introduced infinitesimal analysis in ship building theory. From  a more general point of view, Euler had a systematic approach to the problems of physics, always searching for a differential equation that was at the basis of the problem. For this reason, Euler is considered as one of the main founders of mathematical physics. 
 
  Let us say now a few words on Chebyshev.

There are several short biographies of Chebyshev, yet nothing comparable to the extensive literature that exists on Euler's life and works. We recall here a few elements of Chebyshev's life and works. For more details, we refer the reader to the papers \cite{bio}, \cite{Wassilief}, \cite{You}
and \cite{BJ1999}.

Pafnuty Lvovich Chebyshev was born in 1821. He is sometimes called the Russian Archimedes. His father was a retired Army officer who had contributed to Napoleon's defeat in his attempt to conquer Russia.  Pafnuty Lvovich  received his education at home, essentially from his mother and a cousin. He learned French, which was the language spoken among educated Russians. This became later the language in which he wrote a substantial part of his mathematical  papers. At the age of 11, Chebyshev's parents provided him with private lessons by P. N. Pogorelski, a well-known teacher in mathematics and the author of popular books on elementary mathematics. Chebyshev enrolled at the University of Moscow in 1837, at the Department of Physics and Mathematics. Nikolai Dmetrievich Brashman (1796--1866), who was teaching there pure and applied mathematics, was particularly interested in engineering and hydraulic machines, and Chebyshev was influenced by him.  While he was a student, Chebyshev wrote a paper on iterative methods on approximate solutions of equations. Approximation theory remained one of his favorite subjects until his death.

From the biography contained in his \emph{Collected Works} \cite{T-oeuvres},\footnote{This edition does not contain the complete works of Chebyshev.} we learn that when Chebyshev settled in Saint Petersburg, in 1847, his financial situation was critical. His parents, who used to be wealthy, had lost their properties a few years before, especially during the 1841 famine that affected Russia. They were not able to help him anymore, and his salary as an adjunct-professor was modest. The biographer says that for that reason, Chebyshev became very thrifty, and stayed so throughout his life, but the only thing for which he never spared money was the materials he needed to construct his mechanical models. He was able to spend hundreds of thousands of roubles for his machines.

 Chebyshev became familiar with Euler's writings when he worked on an edition of Euler's papers on number theory, under the direction of Viktor Yakovlevich Bunyakovsky (1804--1889), a project supported by the Academy of Sciences of Saint Petersburg.  Bunyakovsky, during his studies, like several Russian mathematicians of his generation, spent some years in France.\footnote{During a stay in Paris, Bunyakovsky wrote three different doctoral dissertations under Cauchy. The subjects were (1) the rotary motion in a resistant medium of a set of plates of constant thickness and defined contour around an axis inclined with respect to the horizon;  (2) the determination of the radius vector in elliptical motion of planets; (3) the propagation of heat in solids.}   In 1859, Chebyshev was nominated member of the Academy of Sciences  of Saint Petersburg, on the chair of applied mathematics. Bunyakovsky, Mikha\"\i l  Ostrogradskii (1801--1861) and Paul Heinrich Fuss (1798--1855), the great-grandson of Euler, were the three members of the chair of pure mathematics.
 Euler's edition, containing his 99 memoirs on number theory, appeared in 1849 \cite{Euler-C}. Working on Euler's works edition acted certainly as a  motivation for Chebyshev's own research on number theory. In the same year, Chebyshev defended his doctoral thesis in mathematics; the subject was number theory. In 1852, he  published two memoirs in the \emph{Journal de Liouville}, \cite{Cheb-n1} and \cite{Cheb-n2}\footnote{The memoirs were published in that journal in 1852. According to Vassilief  (\cite{Wassilief}  p. 47), they were written in 1848 and 1850 respectively.} on the distribution of prime numbers. 
This was, at that time, one of the most difficult questions in number theory.

Chebyshev's biographer in \cite{bio} says that he thoroughly studied the works of great mathematicians like Euler, Lagrange, Gauss and Abel, and that as a general rule, he avoided reading the works of his contemporary mathematicians, considering that this would be an obstacle for having original ideas. Unlike most of mathematicians of the same epoch, Chebyshev did not like to communicate by letters. On the other hand, he used to travel a lot and the list of mathematicians with which he discussed is impressive. This is the subject of the paper \cite{BJ1989}. According to Vassilief \cite{Wassilief} and Poss\'e \cite{bio}, Chebyshev used to spend almost every summer in Western Europe, especially in France, and particularly in Paris. He published several papers in France, Germany and Scandinavia: at least 27 papers in French journals (among which 17 in Liouville's journal), 3 papers in Crelle's journal, and 5 papers in Acta Mathematica.

   Like Euler, Chebyshev was interested in several fields of mathematics. He published on geometry, algebra, analysis, differential equations, integral calculus, probability, cartography and astronomy, and he had a very strong interest in mechanical engineering. He conceived machines and devices. He was aware of the important place that mathematics holds in the applied sciences, and  conversely, he knew that several practical problems acted  as a motivation for theoretical research. He was especially interested in applied mechanics, and in particular in machine conception, including steam engines, and other engines that transmit motion. The recent survey on linkages, by A. Sossinsky \cite{Sossinsky}, contains some information of Chebyshev's   work on linkages and hinge mechanisms, a subject which he studied thoroughly. He was particularly interested in machines that change a rotational motion into a rectilinear one. See e.g. his paper \cite{Cheb-para}. It is believed that Chebyshev's work on the  approximation of functions was motivated in part by his interest in hinge mechanisms; cf. \cite{You}.
 Chebyshev's \emph{Collected works} \cite{T-oeuvres} contain a report on a 3-month stay he made in France, \cite{voyage}. This stay started on June 21, 1852. At that time,  Chebyshev was adjunct professor at the University of Saint Petersburg, and from there he obtained the imperial permission to spend three months abroad. It appears from the report \cite{voyage} that  Chebyshev  spent much more time in visiting industrial plants and studying machines, than in meeting mathematicians and working on problems of pure mathematics.  
In his report, he describes his observations of 
windmills in Lille and of several machines and models in the \emph{Conservatoire des arts et m\'etiers} in Paris, 
his visit to the
metallurgical plant in Hayange, and to the
paper mills in Coronne and the suburb of Angoul\^eme. 
He also reports on the visit of the
Gouvernment foundry and to the cannon factory in Ruelle,
a visit to a  turbine in a windmaill in Saint-Maur, a
water mill in Meaux,
and an 
arms factory in Ch\^atelleraut. During that visit to France, Chebyshev
 met several eminent French mathematicians, and he spent several evenings discussing with them, only  after his visits to factories in daytime. In this trip, he was able to discuss with Bienaym\'e, Cauchy, Liouville,  Hermite, Lebesgue, Poulignac, Serret and other eminent mathematicians. Chebyshev also made a small trip to London where he discussed with Cayley and Sylvester. He also visited the Royal Polytechnic Institute, where models of various machines were presented. 
We also learn from the same report that during that stay, Chebyshev went to
Brussels where he visited the museum of engines. He was particularly interested in agricultural machines and steam engines.
On his way back to Russia,  Chebyshev made a stop in Berlin and had several discussions with Dirichlet.
A calculating machine is preserved today at the Conservatoire des arts et m\'etiers in Paris. This machine is described in his \emph{Collected works}; see \cite{Cheb-machine}.

   The versatility of Euler and of Chebyshev's mathematical interests are comparable, even if Euler's written corpus is much larger  than that of Chebyshev. 
   Both men were also interested in elementary mathematics, and they wrote reviews and expository papers.

  \section{Cartography}

  We now survey the works of Euler and Chebyshev on Cartography. At the same time, we shall mention works of Lambert, Lagrange and Gauss and a few others on this subject. A detailed history of cartography, in relation with the modern theory of quasiconformal mappings, is contained in the paper \cite{Papa-qc}.
  
At the Saint Petersburg Academy of Sciences, Euler, besides being a mathematician, held the official position of cartographer. He was part of the team of scientists who were in charge of the large-scale project of drawing maps of the new Russian Empire.  Euler also wrote theoretical memoirs  on geography. We mention in particular his \emph{Methodus viri celeberrimi Leonhardi Euleri determinandi gradus meridiani pariter ac paralleli telluris, secundum mensuram a celeb. de Maupertuis cum sociis institutam}
(Method of the celebrated Leonhard Euler for determining a degree of the meridian, as well as of a parallel of the Earth, based on the measurement undertaken by the celebrated de  Maupertuis and his colleagues). This memoir was presented to the Academy  of Sciences of Saint Petersburg in 1741 and published in 1750 \cite{E132}. In 1777, motivated by the practical question of drawing geographical maps, Euler published three memoirs on mappings from the sphere to the Euclidean plane. The memoirs are titled \emph{De repraesentatione superficiei sphaericae super plano} (On the representation of  spherical surfaces on a plane) \cite{Euler-rep-1777},  \emph{De proiectione geographica superficiei sphaericae} (On the geographical projections of spherical surfaces) \cite{Euler-pro-1777}  and  \emph{De proiectione geographica Deslisliana in mappa generali imperii russici usitata} (On Delisle's geographic projection used in the general map of the Russian empire)  \cite{Euler-pro-Desli-1777}. The title of the last memoir  refers to Joseph-Nicolas Delisle (1688-1768), a leading French astronomer and geographer who was also working at the  Academy of Sciences of Saint Petersburg.\footnote{Delisle was a French geographer, and he was part of the team of eminent scientists from several countries of Western Europe who were invited by Peter the Great at the foundation of the  Russian Academy of Sciences. We recall that the monarch  signed the foundational decree of his Academy on February 2, 1724. The group of scientists that were present at the opening ceremony included the mathematicians Nicolaus and Daniel Bernoulli, and Christian Goldbach. Delisle joined the Academy in 1726, that is, two years after its foundation, and one year before Euler's arrival to Saint Petersburg. He was in charge of the observatory of Saint Petersburg. This observatory, situated on the Vasilyevsky Island,  was one of the finest in Europe. During his stay in Russia, Delisle had also access to the most modern astronomical instruments. He was also in charge of drawing maps of the Russian empire. Delisle stayed in  Saint Petersburg from 1726 to 1747. In the first years following his arrival to Russia, Euler assisted Delisle in recording astronomical observations which were used in meridian tables.  In 1747, Delisle left Russia and returned to Paris. He founded there the famous observatory at the h\^otel de Cluny, thanks to a large amount of money he had gathered in Russia. He also published there a certain number of papers containing informations he accumulated during his various voyages inside Russia and in the neighboring territories (China and Japan). Because of that, he had the reputation of having been a spy. In the memoir \cite{Euler-pro-Desli-1777}, Euler calls Delisle ``the most celebrated astronomer and geographer of the time."}

 The three memoirs of Euler on cartography contain important results and techniques of differential geometry.

In the first memoir, \cite{Euler-rep-1777}, Euler examines several projections of the sphere and he systematically searches for the partial differential equations that these mappings satisfy. He starts by recalling that there is no ``perfect" or ``exact" mapping from the sphere onto a plane, and therefore, one has to look for best approximations. He writes:  ``We are led to consider representations which are not similar, so that the spherical figure differs in some manner from its image in the plane." He  highlights the following three kinds of maps:

\begin{enumerate}
\item \label{map1} Maps where the images of all the meridians are perpendicular to a given axis (the ``horizontal" axis in the plane), while all parallels are sent to parallel to it. 

\item \label{map3} Maps which are conformal, that is, angle-preserving.

\item \label{map2} Maps where surface area is represented at its true size.

\end{enumerate}

In the same memoir, Euler gives examples of maps satisfying each of these three properties. The maps which are obtained using various  projections: a projection of the sphere onto a tangent plane, onto a cylinder which is tangent to the equator, etc. He then studies distance and angle distortion under these maps. At the end of his memoir (\S 60), he claims that these investigations have no immediate practical use:\footnote{We are using the translation by George Heine.}
\begin{quote}\small
 In these three Hypotheses is contained everything ordinarily desired from
geographic  as  well  as  hydrographic  maps.   The  second  Hypothesis  treated
above even covers all possible representations.  But on account of the great
generality of the resulting formulae,  it is not easy to elicit from them any
methods of practical use.  Nor, indeed, was the intention of the present work
to go into practical uses, especially since, with the usual projections, these
matters have been explained in detail by others.
\end{quote}

In the same memoir, Euler gives a proof of the fact that there is no mapping from the sphere to the plane which preserves ratios of lengths (\S 9). 

In the second memoir, \cite{Euler-pro-1777}, published in the same year, Euler studies projections that are useful for practical applications. He  writes (\S 20):\footnote{George Heine's translation.}
 \begin{quote}\small
 Moreover, let it be remarked, that this method of projection is extraordinarily appropriate for the practical applications required by Geography,
for it does not strongly distort any region of the earth. It is also important
to note that with this projection, not only are all Meridians and Circles of
Parallel exhibited as circles or as straight lines, but all great circles on the
sphere are expressed as circular arcs or straight lines.  Other hypotheses,
which one might perhaps make concerning the function $\Delta$, will not possess
this straightforward advantage.
\end{quote}

In the third memoir \cite{Euler-pro-Desli-1777}, Euler starts again by presenting the problems to which the stereographic projection leads if it is used in the representation of the Russian Empire. The stereographic projection is a radial projection of the sphere from a point onto a plane tangent to the antipodal point. This projection is extensively used in mathematics. Euler expresses the need for a method of projection, in which the three properties are satisfied:
\begin{enumerate}
\item  All the meridians are sent to straight lines; 
\item degrees of latitudes are preserved;
\item  parallels and meridians meet at right angles in the image. 
\end{enumerate}
He declares that the three conditions cannot be simultaneously satisfied, and the question he addresses becomes that of finding a projection where the deviation regarding the degrees of meridian and parallels is as small as possible, while the two other properties are preserved. Euler then describes a method that Delisle used to draw such a map, and he presents in detail the mathematical theory which underlies this construction.

It is useful to recall here that Euler published several memoirs on spherical geometry, a subject closely related to geography. In fact, theoretical geography uses spherical trigonometry. We also recall that several young collaborators and followers of Euler worked on geography and the drawing of geographical maps. We mention in particular F. T. von Schubert,\footnote{Friedrich Theodor von Schubert (1758-1825) was one of Euler's young collaborators who, after the death of Euler, became the director of the astronomical observatory of the Academy. Schubert was, like Euler, the son of a protestant pastor. His parents, like Euler's parents, first wanted him to study theology and to become a pastor.  Schubert did not follow that path and he studied mathematics and astronomy, without any teacher. He eventually became a specialist in these two fields. He left his native country, Germany, and he became a private mathematics teacher in Sweden, then moved to Estonia. In 1785, two years after Euler's death, he was appointed assistant at the  Academy of Sciences of Saint Petersburg, at the class of geography. In 1789 he  became full member of the Academy.} one of Euler's direct followers, who became a specialist of geography and spherical geometry. We refer the reader to the article \cite{Papa-Inde2} for more information about Schubert. The papers by Schubert on geography include \cite{S0}, \cite{S01} and \cite{S00}. 
Let us also note that some of Euler's memoirs on astronomy concern in fact spherical trigonometry with its practical applications, see e.g. \cite{E14} and \cite{E50}, and there are many others (The Euler Opera Omnia include at least 54 papers on astronomy.).

There are many mathematical and cultural reasons for which the name of Lagrange is associated to the one of Euler. They had several common interests, although they regarded the problems from different points of view. Lagrange was 29 years younger than Euler. Their personalities were quite different, and they had different tastes in mathematics, but they often worked on the same problems. When Lagrange wrote to Euler, on August 12, 1855, to tell him about his new ideas on the calculus of variations, he was only 19 years old. In his approach to the subject, he replaced delicate geometric arguments of Euler by analytic arguments which led directly to what became known later on as the Euler-Lagrange equation.
Euler had a great admiration for his young colleague, who succeeded him as director of the mathematical class at the Berlin Academy of Sciences.

Two years after Euler, Lagrange wrote two important memoirs on cartography \cite{Lagrange1779}. In these memoirs, Lagrange starts by an exposition of the major known geographical maps, in particular those which were known to Ptolemy.  He reviews in particular the stereographic projection, and he highlights the following two properties  it satisfies:

\begin{enumerate}
\item 
 Circles of the sphere are sent to circles in the plane. 
 \item Angles are preserved.
\end{enumerate}

Lagrange then says  that one may consider the geographical maps from a more general point of view, namely, as arbitrary representations of the surface of the globe, and not necessarily obtained by radial projections from points. Euler had a similar point of view.  Lagrange writes  (p. 640) that ``the only thing we have to do is to draw the meridians and the parallels according to a certain rule, and to plot the various places relatively to these lines, as they are on the surface of the Earth with respect to the circles of longitude and latitude." He refers to the work of Lambert, who was probably the first to address the question of characterizing the angle-preserving mappings from the sphere to the plane. Lambert investigated this problem in his \emph{Beitr\"age zum Gebrauche der Mathematik und deren Anwendung} (Contribution to the use of mathematics and its applications) \cite{Lambert-Bey}. Lagrange recalls that Euler, after Lambert, gave a solution of the same problem, and he then gives his own new solution. He considers in detail the case where the images of the meridians and the parallels are circles.

Talking about cartography, one has to say at least a few words about Gauss, who made substantial advances on this subject  using the differential geometry of surfaces, a subject which he developed, motivated precisely by the problem of drawing geographical maps.

  In 1825, Gauss published an important paper, with the title \emph{ Allgemeine Aufl\"osung der Aufgabe: die Theile einer gegebnen Fl\"ache auf einer andern gegebnen Fl\"ache so abzubilden, da\ss \  die Abbildung dem Abgebildeten in den kleinsten Theilen \"ahnlich wird} (General solution of the problem: To map a surface so that the image similar in the smallest parts to what is being mapped) 
   \cite{Gauss-Copenhagen}. The paper concerns the problem of constructing conformal maps. It is known that this paper was influential on Riemann, who was Gauss's student. Gauss  declares in the introduction of that paper that his aim is to construct geographical maps.  He proves several results, among them the fact that every sufficiently small neighborhood  of a point in an arbitrary real-analytic surface can be mapped conformally onto a subset of the plane.  In these investigations, Gauss was motivated by one of his major practical activities, namely, land surveying. Gauss included results of angles and distances that he measured between various points on the Earth  in his famous work \emph{Disquisitiones generales circa superficies curvas} (General investigations on curved surfaces) (1827). It is the task of land surveying that led Gauss eventually to the investigation of triangulations of surfaces,  to the method of least squares (1821), and more generally to the investigation of the differential geometry of surfaces. One should mention in this respect Gauss's result in his \emph{Disquisitiones}  \cite{Gauss-English}, which he calls the ``remarkable theorem," 
saying that the parameter known today as Gaussian curvature is the obstruction for a sphere to be faithfully represented on the plane \S 12 (cf. p. 20 of the English translation).

It was natural that Chebyshev, who was, among other things, a geometer, an analyst and an applied mathematician, and who was, as we already noted, a devoted reader of Gauss, became interested in cartography.

  Chebyshev's \emph{Collected papers} edition \cite{T-oeuvres} contains two papers on cartography, \cite{Cheb1} and \cite{Cheb2}.  The first paper is more technical (it contains several formulae), and in some sense it is a sequel to the work of Lagrange on the same subject. The second paper, which is longer, contains more theoretical and philosophical considerations, as well as some historical notes on the discovery of differential calculus. On p. 14 of his biography \cite{Wassilief}, Vassilief  says that it is probably in reading Euler that Chebyshev became interested in the drawing of geographical maps.

In the introduction of the first paper, \cite{Cheb1}, Chebyshev recalls that it is easy to construct geographical maps that preserve angles (he says: ``such that the infinitesimal elements on the sphere and their representation are preserved"), but that for these maps, length is distorted, and the ratios of length elements between points of the sphere and its representation vary from point to point. Thus, one may conceive  that there exist maps where these deviations are the smallest possible. This is the problem in which he is interested.

  Chebyshev says that from the mathematical point of view, this question bears some strong analogy with some problems he considered before, concerning linkages. He calls the study of linkages the \emph{theory of mechanisms, known under the name of parallelograms}. The word \emph{parallelogram} refers here to four-bar linkages that have the form of parallelograms that appear in the machines studied by Chebyshev. Let us note that this word appears several times in the present paper, in various contexts. 
  
  In the second paper \cite{Cheb2}, Chebyshev returns again to the analogy, and he talks about Watt's parallelograms, a mechanical linkage described in 1784 by James Watt (1736--1819) when he patented the so-called \emph{Watt steam engine}.
Chebyshev mentions the fact that in both theories (geographical maps and linkages), one looks for a function of two variables that realizes a minimum among functions which satisfy a certain partial differential equation. Thus, he includes the two theories in the setting of the calculus of variations, a subject which was dear to Euler and Lagrange.

Chebyshev bases his investigations on a formula  found by Lagrange, for a quantity he calls the \emph{magnification ratio}. This is the ratio between a length element at a point on the sphere and its image by a map. The setting is the one of Lagrange in his paper on cartography \cite{Lagrange1779}, in which he studies the maps of the sphere into the Euclidean plane that are "similarities at the infinitesimal level." These are the angle-preserving mappings that were thoroughly studied by Euler, and then by Lagrange. In particular, for such mappings, the magnification ratio, depends only on the point and not on the chosen direction.
Lagrange gave the following formula for the magnification ratio, which Chebyshev recalls:
\[m=\frac{\sqrt{f'(u+t\sqrt{-1})F'(u-t\sqrt{-1})}}{\frac{2}{e^u+e^{-u}}}.
\] 
Here, $f$ and $F$ are arbitrary functions.
The formula gives 
\[\log m = \frac{1}{2}\log \big( f'(u+t\sqrt{-1})\big)+  \frac{1}{2}\log \big( F'(u-t\sqrt{-1})\big) -\log \frac{2}{e^u+e^{-u}}.\]

The first two terms of the right hand side of this equation, involving arbitrary functions, constitute the solution of a partial differential equation of the form
\[\frac{\partial^2 U}{du^2}+\frac{\partial^2 U}{dt^2}=0
\]
which is nothing else than the Laplace equation. 
In the rest of his paper, Chebyshev solves the question of finding the functions that have the smallest magnification ratio.
He establishes a relation between a quantity which Lagrange calls the ``projection exponent" and the shape of the curve that bounds the country which the map is meant to represent.

 The second paper is a sequel to the first. Using the computations in the first paper, Chebyshev  
  \cite{Cheb2}, determines, for a given country, the center of the projection and the value of this exponent which are best suited for the drawing of a geographical map. 
  
  At the beginning of the second paper  \cite{Cheb2}, Chebyshev talks about the mutual influence of theoretical and practical sciences. He declares that every new field of study in theoretical mathematics originates in practical  problems, and that a particularly important class of problems arises from questions requiring the maximization or the minimization of certain quantities. For instance, this is how differential calculus, and later, the calculus of variations, were born. He declares that the more practical needs are demanded, the further theory goes. He gives examples from Newton's \emph{Principia}. The construction of geographical maps is another very good example of this link between practical and theoretical sciences. He considers, as in the first paper, the magnification ratio, a function defined on the region for which one wants to draw a map.  He obtains the following theorem (\cite{Cheb2} p. 242):
  \begin{quote}
   \emph{The best projection of a country, or a region, in the sense of the smallest length distortion consists of one for which the length distortion ratio on the boundary of the region to be represented is constant}.
   \end{quote}
   Chebyshev also provides a method for the computation of this distortion ratio.

  With this theorem in mind, Chebyshev notes that the search for the best map amounts to the solution of a certain partial differential equation defined on a region with given boundary values. He mentions that the theory  developed from the partial differential equation point of view is analogous of the theory of heat propagation. Indeed, once the value of the 
   magnification ratio on the boundary of the country is known, the rest is found by solving Laplace's equation. 
   
   A proof of Chebyshev's theorem, based on the same idea (the solution of the Laplace equation),  is contained in Milnor's paper \cite{Milnor}, p. 1111. Before Milnor, Darboux gave a solution of the same problem, also based on Chebyshev's ideas, using potential theory, \cite{Darboux-Chebyshev}. It seems that Milnor was not aware of the work of Darboux on the question. 
   
Chebyshev's paper \cite{Cheb2} ends with some explicit examples of maps.

 Some of the works started by Chebyshev on the drawing of geographical maps were continued by his student  Dmitry Aleksandrovich Grav\'e (1863--1939), who attended  Chebyshev's courses in the early 1880s. Grav\'e also translated from Russian into French Chebyshev's papers on geographical maps 
\cite{Cheb1} and \cite{Cheb2} that are included in the French edition of his \emph{Collected works} \cite{T-oeuvres}. The subject of Grav\'e's doctoral dissertation was \emph{On the main problems of the mathematical theory of construction of geographical maps}. He defended it in Saint Petersburg in 1896 (this was after Chebyshev's death). This work concerns equal area projections of the sphere, and it is based on ideas of Euler, Lagrange and Chebyshev, cf. \cite{Erm} and \cite{Mes}. Concerning Chebyshev's courses, Grav\'e  wrote (cf. \cite{Bernstein}): 
 \begin{quote}\small 
 Chebyshev was a wonderful lecturer. His courses were very short. As soon as the bell sounded, he immediately dropped the chalk, and, limping, left the auditorium. On the other hand he was always punctual and not late for classes. Particularly interesting were his digressions when he told us about what he had spoken outside the country or about the response of Hermite or others. Then the whole auditorium strained not to miss a word.
 \end{quote}
Let us also note another account of Chebyshev's courses, by Aleksander Mikhailovich Lyapunov (1857--1918), the mathematician well known for his work on the stability of dynamical systems. Lyapunov attended Chebyshev's lectures in the 1870s. The quote is again from \cite{Bernstein}: 
 \begin{quote}\small 
 [Chebyshev's]  courses were not voluminous, and he did not consider the quantity of knowledge delivered; rather, he aspired to elucidate some of the most important aspects of the problems he spoke on. These were lively, absorbing lectures; curious remarks on the significance and importance of certain problems and scientific methods were always abundant. Sometimes he made a remark in passing, in connection with some concrete case they had considered, but those who attended always kept it in mind. Consequently his lectures were highly stimulating; students received something new and essential at each lecture; he taught broader views and unusual standpoints.
  \end{quote}

Let us finally mention, considering the question of geographical maps and projections of a curved surface onto a Euclidean plane, that  Beltrami studied thoroughly the general question of the characterization of the mappings from a given surface onto the plane that send geodesics to Euclidean straight lines.
This problem is closely related to Hilbert's Problem IV, which he presented among the list of 23 problems he compiled on the occasion of his talk at the Paris 1900 ICM \cite{Hilbert-Problems}, in which he asks for the study of the metrics on the plane whose geodesics are the Euclidean straight lines.

  \section{Chebyshev on the fitting of garments}
  
 We now consider a second problem that Chebyshev studied, namely, the problem of the fitting of garments. Chebyshev presented this problem at a meeting of the Association Fran\c caise pour l'Avancement des Sciences.\footnote{The \emph{Association Fran\c caise pour l'Avancement des Sciences} (French Association for the Advances of Sciences) was founded in 1872. Its initial goal was to make links between researchers in different fields, and also between scientists and science amateurs, coming  from different backgrounds. Its aims also included the popularization of science. The presidents of this society were elected every year, and they came from various fields of sciences. They included in particular the mathematicians Paul Appell (1908),  Emile Borel (1925), Elie Cartan (1933) and Paul Montel (1946). Between the years 1873 and 1882, Chebyshev presented sixteen reports at various sessions of the Association. Among them are his report \emph{Sur la coupe des vêtements} (On the cutting of garments)  which we consider here and the one on the calculating machine which we already mentioned and which is described in his paper \cite{Cheb-machine}.} In his talk, Chebyshev starts by declaring that the idea of this work originates from a communication made at the same Association two years earlier by \'Edouard Lucas, on the geometry of weaving fabrics with rectilinear wires.\footnote{The mathematician \'Edouard Lucas (1842--1891) was 21 years younger than Chebyshev. The two men attended the 1876 session of the Association Fran\c caise pour l'Avancement des Sciences which took place in Clermont-Ferrand. Lucas is known for his work on number theory. His name is associated to  \emph{Lucas numbers}, which are related to Fibonacci numbers by a recurrence relation.  It is conjectured that there are infinitely many prime Lucas numbers. Lucas is also known for his work on satin squares, which makes connections between number theory and fabrics. His writings in this domain have both a mathematical and general public character, cf. \cite{Lucas0} \cite{Lucas1} \cite{Lucas2}. For an exposition of Lucas' work on fabrics, in relation with other nineteenth-century works by various people, we refer the reader to the paper \cite{Decaillot} by Decaillot. Lucas died at the age of 49, after an accident at the banquet of the annual congress of the Association Fran\c caise pour l'Avancement des Sciences. A waiter  dropped a pile of plates, and a piece of broken plate cut Lucas on the cheek. This caused a skin inflammation and Lucas died a few days later. Let us note that the \emph{Gauss-Lucas theorem}, saying that the convex hull in the complex plane of the roots of a polynomial contains the roots of the derivative polynomial,  refers to F\'elix Lucas, a contemporary of \'Edouard Lucas. We mention by the way a geometric proof of that theorem due to William Thurston, published recently, see \cite{Cheritat-Lei}. We shall talk about Thurston's work, in relation with that of Chebyshev, in the last section of this paper.}    In fact, Lucas thoroughly studied questions related to garments, but it is not clear that he considered the specific one in which Chebyshev was interested. The paper \cite{Lucas2} by Lucas, which was published in 1880, is based on a previous paper he wrote. It concerns weaving styles,  more particularly the various ways in which the two perpendicular threads that constitute a piece of fabric intersect and give rise to different kinds of textile. In that paper, Lucas declares that his work on the subject gave rise to works by several people and he mentions among them the memoir \cite{Chebyshev} by Chebyshev. He considers these works as new applications of number theory. It is not clear that the paper of Chebyshev that we are analyzing here falls into this category, and the relation with number theory is not visible.
     Even though the subject matter of Chebyshev's paper is related to garments, the content is very different from the subject in which Lucas was interested in. The point we want to make here is that even though Chebyshev declares that the idea of his investigation occurred to him after he heard the talk by Lucas, there is no relation between the problems considered by the two mathematicians, except for the fact that they involve fabrics.
     
The problem addressed by Chebyshev is the following. One starts with a fabric made out of a net of two perpendicular threads. In the language of fabric industry, these two threads are called warp threads and weft threads. Each two perpendicular threads of the net are attached at their intersection point. When such a piece of fabric is used to cover a part of a human body -- for instance the head, in the case of head nets -- it takes the form of that part of the body. The threads are rigid at their intersection points, that is, the lengths between any two consecutive intersection points is fixed. The flexibility is only at the level of the angles at these intersections. While at the beginning, the treads make Euclidean rectangles,  after the deformation, the rectangles become curvilinear parallelograms. Thus, we can can see the relation with  the theory of linkages, which, as we already mentioned, Chebyshev called the \emph{theory of mechanisms, known under the name of parallelograms}. The question which Chebyshev addresses, in the context of woven fabric, is to see whether with such an amount of flexibility the piece of fabric can take the exact form of the part of body to which it is designated. Another question, which he considers later on, is to study the optimal curves and their form, along which several pieces of fabric might be sewed, in order to cover the given part of the body.

 In his communication to the Association, Chebyshev described the problem, and solved it for the case where the surface that has to be covered is the hemisphere. The paper includes several practical examples. Chebyshev gave formulae for the solution and, with his usual faithfulness to one of his favorite subjects, approximation theory, he explained how one can obtain approximate solutions.\footnote{Tikhomirov \cite{Ti}  considers that ``[Chebyshev] set the foundations of the Russian school of approximation theory."}  He also studied the  question of sewing several pieces of fabric, that is, he replaced the problem by the more general one that consists of binding a surface with several pieces of fabric sewed together. Chebyshev found a differential equation satisfied by the Gaussian curvature of a surface covered by the cloth. 
The term \emph{Chebyshev net} was given to such a net, that is, a net made out of quadrilaterals of fixed side length. The question becomes that of the possibility of covering an arbitrary surface with a Chebyshev net. The local problem was solved affirmatively by Bieberbach  in 1926 \cite{Bieberbach}. It is known now that for the global problem, one needs to deal with singularities. Burago, Ivanov and Malev proved in \cite{BIM} a general result on the existence of a covering by a Chebyshev net for complete simply connected Alexandrov surfaces under some constraints on the curvature. 

In his paper \cite{Chebyshev}, Chebyshev included the problem of the fitting of garments in the setting of the differential geometry of surfaces. From this point of view, a Chebyshev net on a surface provides a (local) coordinate system in which the horizontal and the vertical vector fields have norm one, and the coordinate local transformations are translations. 
Chebyshev noted that in the coordinates associated to the curves forming the net, the line element on the surface can be written in the form
 \begin{equation} \label{eq:T} 
 ds^2= du^2+dv^2+2\cos\varphi dudv,
\end{equation}
 where $u$ and $v$ are the parameters along the curves of the net and $\varphi$ the angle between these curves. 
  The equation satisfied by the Gaussian curvature $K$ of a surface covered by a  Chebyshev net is then of the form
\[ \varphi_{uv}+K\sin \varphi =0.
\] 

It turns out that for surfaces of constant curvature, this is precisely a Sine-Gordon equation; cf. \cite{Rosenfeld}.  The authors of \cite{Rosenfeld} note that
 Hilbert proved in \cite{H1900} that the Chebyshev net is formed by its asymptotic curves, that is, the curves that are tangent at each point to the \emph{principal tangents} (the tangent to the directions where the magnification ratio is maximal and minimal). This shows that the angle between the asymptotic curves on such a surface satisfies a Sine-Gordon equation. 
 
 There is a large number of relatively recent papers studying Chebyshev nets.  There are also applications  in industry. 
 
 The French original  version of Chebyshev's \emph{Collected works} \cite{Chebyshev}  contains only a summary of Chebyshev's original talk.  The complete text appears in the Russian version of the Collected Works,  \cite{Chebyshev-R},  vol. 5, p. 165-170. A French version of this extended text is included in the PhD thesis \cite{Decaillot-T} of Decaillot and in the paper \cite{Ghys} by Ghys. The last paper is a nice update of Chebyshev's  problem of the fitting of garments.

\begin{figure}[htbp]
\centering
\includegraphics[width=6.5cm]{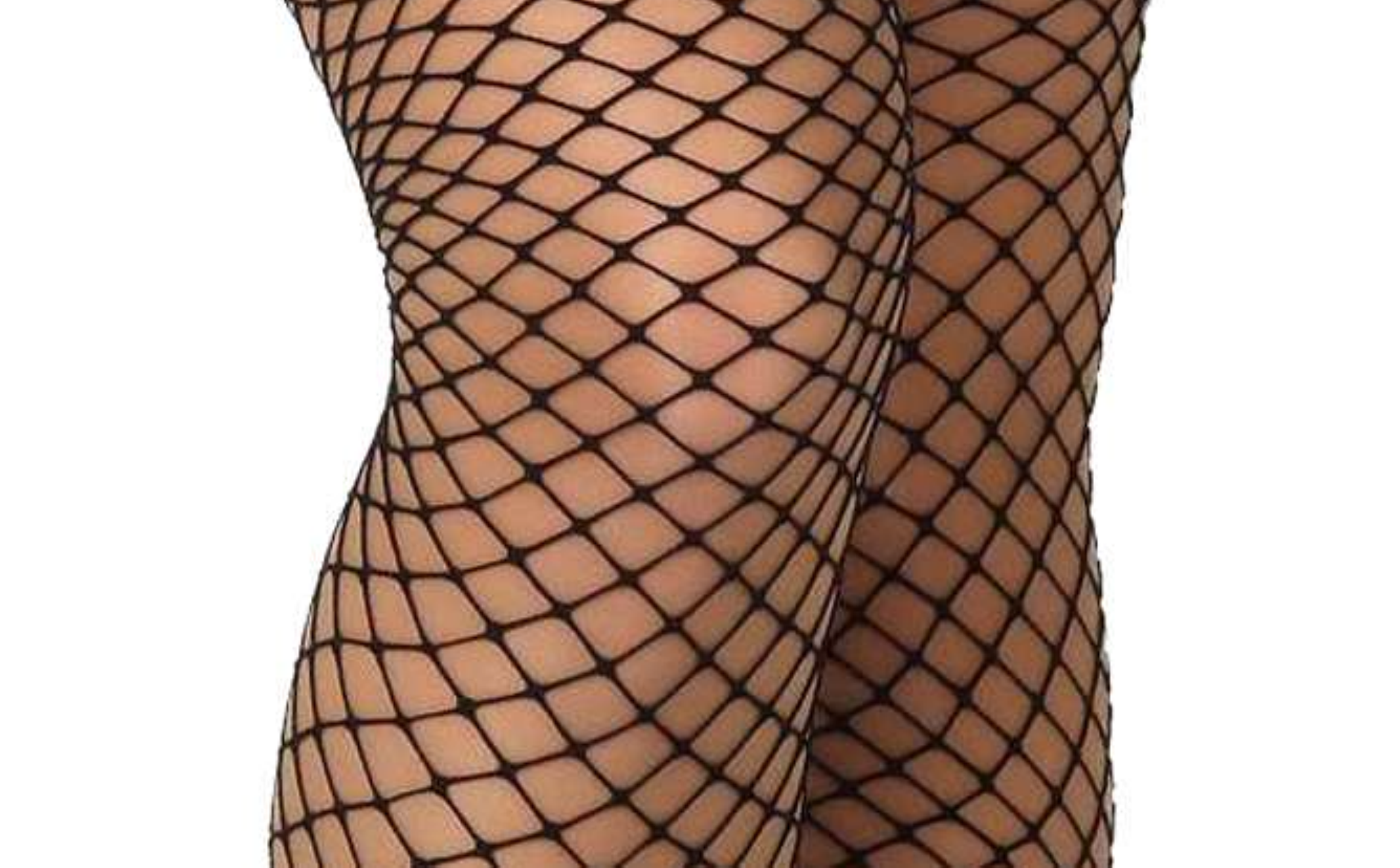}
\caption{A Chebyshev net takes the form of the surface on which it is applied.} \label{pp}
\end{figure}

 We mentioned in the introduction that the problem of drawing geographical maps and Chebyshev's problem of the fitting of garments are in some sense inverse to each other. The comparison between these problems may be pushed further. On the sphere, one has a natural system of orthogonal coordinates: the equator together with its parallels, which form a foliation (with two singular points) of the sphere, and the meridian lines, which form an orthogonal (singular) foliation. On a geographical map, one would like to represent the image of these two orthogonal foliations; the result is a pair of transverse foliations of (a subset of) the Euclidean plane which, in the general case, are not necessarily orthogonal.
 In the problem of the fitting of garments, one starts conversely with two orthogonal foliations on the Euclidean piece of fabric, the warp threads and the weft threads, and one is interested in the images of these foliations (which are no longer orthogonal) on the curved surface that is covered.  In Figure \ref{pp}, an example of Chebyshev net is given. To stay in the setting of Chebyshev's theory, one has to make the assumption that the thread used for the fishnet is not elastic.
 
 Darboux, in his ICM talk in Rome (1908) \cite{Darboux-ICM}, mentions questions related to geographical maps and he also mentions the work of Chebyshev on fitting garments. He includes both subjects in his paper, titled \emph{Les origines, les m\'ethodes et les probl\`emes de la g\'eom\'etrie infinit\'esimale} (The origins, methods and problems of infinitesimal geometry). He presents the following question, which is slightly more general than the one studied by Chebyshev. Consider a piece of fabric, in the form of a net which ladies use to cover their hair. The net is formed by the two perpendicular families of threads (the warp threads and the weft threads). These threads are attached at their intersections, in such a way that they form small rectangles. The net can be deformed in such a way that the angles at the vertices of the rectangles may vary, but not the side lengths.  
The question is to determine the form of the net when it is placed on a surface (part of a human body).  The problem is slightly more general  than the one Chebyshev studied. The latter studied the case where the rectangles are squares. 

 Darboux writes the length element of the surface in terms of the angle variable $\theta$ in the following form:
\begin{equation} \label{eq:Darboux} ds^2= A^2d\alpha^2+C^2d\beta^2+2AC\cos\theta d\alpha d\beta,
\end{equation}
 where $A$ and $C$ are known functions of $\alpha$ and $\beta$ and $\theta$ is the unknown. 
 Equation (\ref{eq:T}) which Chebyshev used is a special case of  Equation (\ref{eq:Darboux}), obtained by taking $A=C=1$. Darboux discusses this question and he says that Chebyshev was satisfied with approximate solutions of the problem. He claims that the latter failed to notice that the general and complete solution of his problem, in the special case of the sphere, is related to the question of determining the surfaces of constant negative curvature.

\section{Euler again}

Finally, we mention a problem which Euler studied that involves the geometry of fabrics. This is the question of sails of ships. Euler first investigated this question in 1727, at the age of 20, in a paper he sent to the French Academy of Sciences, as a solution to a problem which the Academy proposed as a contest \cite{Euler-Meditationes}.  The title of Euler's essay is \emph{Meditationes super problemate nautico, quod illustrissima regia Parisiensis Academia scientarum proposuit} (Thoughts on a nautical problem, proposed by the illustrious Royal Academy of Sciences in Paris). 
 The problem proposed by the Academy asks for the best way to place masts on vessels. The height of the masts, and the height and the width of the sails are to be determined. The study of the problem involves the investigation of the capacity of the sails, on which depend the forces of wind on them. The inclination of the ship also depends on the force of the wind on the sails. Johann Bernoulli,\footnote{This is Joahann I Bernoulli (1667--1748).} who was Euler's teacher at the University of Basel and who encouraged him to work on that problem, had already published a paper on the same subject, in 1714, which he called \emph{Th\'eorie de la man\oe uvre des vaisseaux} (Theory of maneuver of ships). In that paper, Bernoulli  determines the shape of a sail in terms of the pressure exerted by the wind on it.
In the paper submitted to the Academy, Euler extended ideas of Archimedes contained in his treatise \emph{On floating bodies} \cite{Archimedes} and introduced in the subject techniques of differential calculus, in particular partial differential equations. His paper earned a honorable mention. Euler remained interested in these questions for the rest of his life.
Several years later, he completed his major opus on ship building, the \emph{Scientia navalis} (Naval science). This two-volume treatise appeared in 1749 in Saint Petersburg \cite{Euler-Navalis}.

\medskip

  \section{By way of conclusion}

  William Thurston (1946--2012), who can reasonably be considered as the greatest geometer of recent times, revived two subjects  which were dear to Chebyshev, namely,  mechanical linkages and clothing design.

  Thurston did not write up his ideas on linkages, but several people remember his lectures. The subject is mentioned in his paper \cite{TW}, with the expanded French version \cite{TWF}. It is also mentioned in his book \cite{TP}.  In a recent correspondence with the author of the present article, Bill Abikoff writes: ``Thurston was characteristically terse in his discussion of spaces formed by flexible linkages. He either told me directly or someone told me the following. His response to the question of which topological spaces appear as the configuration space of a flexible linkage was: \emph{all}."\footnote{In another recent mail, Abikoff mentions that his father worked with every aspect of the creation of machine-woven fabric in most of his professional life. Another thing he mentions is that a couple of years ago he had a discussion with Louis Nirenberg about the contemporary relevance of Euler's work on sails of ships in relation with compressible fluid flow. They both felt that there was renewed interest in this topic because of wind farms: another relation with Chebyshev's work!} One may note that this is a vast generalization of a result of Kempe which says that any bounded piece of an algebraic curve is drawable by some linkage, cf.  \cite{Kempe}. One may recall by the way that the definition of a curve, in the period of Greek mathematics, was mechanical (i.e. in terms of linkages). Thus, through equations and Descartes' geometry, one does not find more curves than those which were considered in Greek antiquity.

  Regarding Thurston's contribution to linkages, see also the exposition in \cite{Sossinsky}.

   In his important paper on minimal stretch maps between hyperbolic surfaces (1985), Thurston alludes to the relation between cloth deformation and some extremal maps between surfaces. After he states the main problem he studies in this paper -- the problem of finding best Lipschitz maps between surfaces, he writes: 
   \begin{quote}\small This is closely related to the canonical problem that arises when a person on the standard American diet digs into his or her wardrobe of a few years earlier. The difference is that in the wardrobe problem, one does not really care to know the value of the best Lipschitz constant -- one is mainly
concerned that the Lipschitz constant not be significantly greater than 1. We
shall see that,  just as cloth which is stretched tight develops stress wrinkles, the least Lipschitz constant for a homeomorphism between two surfaces is dictated by a certain geodesic lamination which is maximally stretched...
\end{quote} 
 Thurston elaborated on the comparison between extremal maps between surfaces and deformation of clothes in several of his lectures.

  Around the year 2010, Thurston started working, together with the Japanese fashion label Issey Miyake,  on the conception of geometrically-inspired garments, in particular, patterns based on non-Euclidean geometry. In a well-known interview (ABC News report) titled \emph{Fashion and advanced mathematics meet at Miyake: Fashion and advanced mathematics collide at Japanese label Issey Miyake}, along with the  Japanese designer Dai Fujiwara, Thurston says: ``I have long been fascinated (from a distance) by the art of clothing design and its connections to mathematics." In 2010,  Miyake presented in Paris a fashion collection with the title ``8 Geometry Link Models as Metaphor of the Universe," inspired by Thurston's work. The collection is also known under the name ``Poincar\'e Odyssey."
 The following is an excerpt of a text of Thurston which was available at the fashion show:
 \begin{quote}\small
 Many people think of mathematics as austere and self-contained. To the contrary, mathematics is a very rich and very human subject, an art that enables us to see and understand deep interconnections in the world. The best mathematics uses the whole mind, embraces human sensibility, and is not at all limited to the small portion of our brains that calculates and manipulates with symbols. Through pursuing beauty we find truth, and where we find truth, we discover incredible beauty.
 \end{quote}  
   
   \medskip
   
   \noindent \emph{Acknowledgements.---} I would like to thank Valerii Galkin, from Surgut State University, who invited me to give a talk at the conference dedicated to Chebyshev (May 17-19, 2016). I would also like to thank Norbert A'Campo,  Bill Abikoff, Vincent Alberge, Jacques Franchi, Elena Frenkel,  Etienne Ghys, Olivier Guichard, Marie-Pascale Hautefeuille, Misha Katz, Ken'ichi Ohshika and Sumio Yamada. All of them read a preliminary version of this paper and made corrections and suggestions.

\printindex

\begin{thebibliography}{00}

\bibitem{Archimedes} Archimedes, On floating bodies. In: The works of Archimedes, ed. T. L. Heath, Reprint, Dover, Mineola, NY, 2002.

 \bibitem{Bernstein} S. N. Bernstein, Chebyshev's influence on the development of mathematics (Russian), Uch. Zap. Mosk. Gos. Univ. 91 (1947), 35--45. English translation by  O. Sheynin in Math. Sci. 26 (2) (2001), 63--73. 

\bibitem{Bieberbach} L. Bieberbach, \"Uber Thebychefsche Netze auf Flächen negativer Krümmung, sowie auf einigen weiteren Flächenarten. Preuss. Akad. Wiss., Phys. math. Kl. 23, 294--321.
\bibitem{BIM}
Y. D. Burago, S. V. Ivanov, and S. G. Malev. Remarks on Chebyshev coordinates, Zap. Nauchn. Sem. S.-Peterburg. Otdel. Mat. Inst. Steklov. (POMI), 329 (Geom. i Topol. 9) 195 (2005) 5--13, 

\bibitem{BJ1989} P. L. Butzer and F. Jongmans, P. L. Chebyshev (1821--1894) and his contacts with Western European scientists, Historia Mathematica 16 (1989) 46--68.


\bibitem{BJ1999} P. L. Butzer and F. Jongmans, P. L. Chebyshev (1821--1894): A guide of life and works. Journal of Approximation Theory, 96 (1999) 111--138.
 
 \bibitem{Cheb-n1} P. L. Chebyshev, Sur la fonction qui d\'etermine la totalit\'e des nombres premiers inf\'erieurs \`a une limite donn\'ee. M\'emoires pr\'esent\'es \`a l'Acad\'emie Imp\'eriale des sciences de Saint-P\'etersbourg par divers savants, t. VI, 1851, p. 141--157 and Journal de math\'ematiques pures et appliqu\'ees, t. 17, 1852, p. 341--365. Reprinted in P. L. Tchebycheff, \OE uvres \cite{T-oeuvres}, Vol. 1,  p. 29--48. Reprint, Chelsea, NY. 
 
 
 \bibitem{Cheb-n2}  P. L. Chebyshev, M\'emoire sur les nombres premiers, M\'emoires pr\'esent\'es \`a l'Acad\'emie Imp\'eriale des sciences de Saint-P\'etersbourg par divers savants, t. VII, 1854, p. 17-33 and Journal de math\'ematiques pures et appliqu\'ees, t. 17, 1852, p. 366--390. Reprinted in P. L. Tchebycheff, \OE uvres \cite{T-oeuvres}, Vol. 1,  p. 51--70. Reprint, Chelsea, NY. 

  

  \bibitem{Cheb-para} P. L. Chebyshev, Th\'eorie des m\'ecanismes connus sous le nom de parall\'elogrammes. M\'emoires pr\'esent\'es \`a l'Acad\'emie Imp\'eriale des sciences de Saint-P\'etersbourg par divers savants, VII, 1854, p. 539--568. Reprinted in P. L. Tchebycheff, \OE uvres \cite{T-oeuvres}, Vol. 1,  p. 111--143, Reprint, Chelsea, NY. 
  
  
\bibitem{voyage} P. L. Chebyshev, Rapport du professeur extraordinaire de l'universit\'e de Saint-P\'etersbourg Tchebychef sur son voyage \`a l'\'etranger. Reprinted in P. L. Tchebycheff, \OE uvres \cite{T-oeuvres}, 
  Vol. 2,  p. vii-xiii.


  
  \bibitem{Chebyshev} P. L. Chebyshev, Sur la coupe des v\^etements, Assoc. Fran\c c. pour l'Avancement des Sciences, 7\`eme session \`a  Paris, 28 Ao\^ou 1878, 154--155. Reprinted in:  P. L. Tchebycheff, \OE uvres \cite{T-oeuvres},  Vol. 2, p. 708 (excerpt). Reprint,  Chelsea, NY. The complete text in Russian  appears in the Russian \emph{Collected Works} \cite{Chebyshev-R}.
  
   \bibitem{Chebyshev-R}  P. L. Chebyshev, Complete collected works (Russian), Acad. of Sciences of the USSR, Moscow-Leningrad, 1944/1951.


  \bibitem{Cheb1} P. L. Chebyshev, Sur la construction des cartes g\'eographiques. Bulletin de la classe physico-math\'ematique de l'Acad\'emie Imp\'eriale des Sciences de Saint-P\'etersbourg, Tome VIV, 1856, p. 257--261. Reprinted in P. L. Tchebycheff, \OE uvres \cite{T-oeuvres} Vol. 1,  p. 233--236, Reprint, Chelsea, NY. 
  
    \bibitem{Cheb2} P. L. Chebyshev, Sur la construction des cartes g\'eographiques. Discours prononc\'e le 8 f\'evrier 1856 dans la s\'eance solennelle de l'Universit\'e Imp\'eriale de Saint-P\'etersbourg (traduit par A. Grav\'e). Reprinted in P. L. Tchebycheff, \OE uvres \cite{T-oeuvres}, Vol. 1, p. 239--247,  Reprint,  Chelsea, NY. 
    
    
        \bibitem{Cheb-machine}  P. L. Chebyshev, Une machine arithm\'etique \`a mouvement continu, Reprinted in P. L. Tchebycheff, \OE uvres \cite{T-oeuvres}, Vol. 2,  p. 237--247, Reprint, Chelsea, NY. 
    
    
    \bibitem{T-oeuvres} P. L. Chebyshev, \OE uvres, edited by A. Markoff and N. Sonin, Imprimerie de l'Acad\'emie Imp\'eriale des Sciences,  Saint Petersburg, 2 volumes, 1899-1907.
    
\bibitem{Cheritat-Lei} A. Ch\'eritat, Y. Gao, Y. Ou and T. Lei, A refinement of the Gauss-Lucas theorem, C. R. Acad. Sci. Paris, Ser. I 315 (2015) 711--715.



\bibitem{Darboux-ICM} G. Darboux, Les origines, les m\'ethodes et les probl\`emes de la g\'eom\'etrie infinit\'esimale. Proceedings of the IVth ICM, Rome. Vol. 1, 105--122 (1909).

\bibitem{Darboux-Chebyshev} G. Darboux, Sur la construction des cartes géographiques, Bulletin des sciences mathématiques,  Bulletin des sciences mathématiques,  35 (1911) 23--28.




\bibitem{Decaillot} A.-M. Decaillot, G\'eom\'etrie des tissus. Mosa\"\i ques. \'Echiquiers. Math\'ematiques curieuses et utiles.
Revue d'histoire des math\'ematiques, 8 (2002), p. 145--206.

\bibitem{Decaillot-T} A.-M. Decaillot, \'Edouard Lucas (1842--1891): le parcours original d'un scientifique fran\c cais dans la deuxième moiti\'e du XIXe siècle, thèse, Universit\'e de Paris 5--Ren\'e Descartes, 2 vol., 1999.
\bibitem{Euler-Meditationes} L. Euler, Meditationes super problemate nautico, quod illustrissima regia Parisiensis Academia scientarum proposuit, Pièce qui a remport\'e le prix de l'Academie Royale des sciences, 1727, 1728, p. 1--48
Opera Omnia, Series 2, Volume 20, p. 1--35.


 \bibitem{E17} L. Euler,  Einleitung zur Rechen-Kunst zum Gebrauch des Gymnasii bey der Kayserlichen Academie der Wissenschafften in St. Petersburg. Gedruckt in der Academischen Buchdruckerey 1738. Opera Omnia,  Series 3, Volume 2, p. 1--304.
 

\bibitem{E14} L. Euler,  Solutio problematis astronomici ex datis tribus stellae fixae altitudinibus et temporum differentiis invenire elevationem poli et declinationem stellae. Auct. L. Eulero, Commentarii academiae scientiarum Petropolitanae 4, 1735, p. 98-101. Opera Omnia, Series 2, Vol. 30, p. 1--4.


\bibitem{E50} L. Euler,  Methodus computandi aequationem meridiei. Auctore L. Eulero 
 Originally published in Commentarii academiae scientiarum Petropolitanae 8, 1741, p. 48--65. Opera Omnia, Series 2, Vol. 30, p. 13--25.


 
\bibitem{E77} L. Euler, New principles of gunnery, London, 1742. German edition, Neue Grundsätze der Artillerie, 1745. Opera Omnia, Series 2, Volume 14, p. 1--409.

\bibitem{E118} L. Euler, Sur la perfection des verres objectifs des lunettes, M\'emoires de l'Acad\'emie des sciences de Berlin 3, 1749, p. 274--296. Opera Omnia, Series 3, Volume 6, p. 1--21 

\bibitem{E132} L. Euler,  Methodus viri celeberrimi Leonhardi Euleri determinandi gradus meridiani pariter ac paralleli telluris, secundum mensuram a celeb. de Maupertuis cum sociis institutam, Commentarii academiae scientiarum Petropolitanae 12, 1750, p. 224--231. Opera Omnia, Series 2, Volume 30, p. 73--88.


\bibitem{E196} L. Euler,  Emendatio laternae magicae ac microscopii solaris, Novi Commentarii academiae scientiarum Petropolitanae 3, 1753, p. 363--380.
Opera Omnia, Series 3, Volume 6, p. 22--37. 

 
 \bibitem{E194} L. Euler,    De machinis in genere
 Novi Commentarii academiae scientiarum Petropolitanae 3, 1753, p. 254--285. 
    Opera Omnia, Series 2, Volume 17, p. 40--65.
    
      \bibitem{E235} L. Euler, Sur l'action des scies
M\'emoires de l'acad\'emie des sciences de Berlin  12, 1758, p. 267--291.
    Opera Omnia, Series 2, Volume 17, p. 66--88.
    
    
          \bibitem{E238} L. Euler, Sur la force des colonnes
M\'emoires de l'acad\'emie des sciences de Berlin 13, 1759, p. 252--282.
    Opera Omnia, Series 2, Volume 17, p. 89--118.
    
              \bibitem{E233} L. Euler, Recherches plus exactes sur l'effet des moulins \`a vent, M\'emoires de l'acad\'emie des sciences de Berlin 12, 1758, p. 165--234.
    Opera Omnia, Series 2, Volume 16, p. 65--125.

\bibitem{E239} L. Euler,   Règles g\'en\'erales pour la construction des T\'el\'escopes et des Microscopes, de quelque nombre de verres qu'ils soient compos\'es.  M\'emoires de l'Acad\'emie des sciences de Berlin 13, 1759, p. 283--322.
Opera Omnia, Series 3, Volume 6, p. 44--73.


 \bibitem{E426} L. Euler, Th\'eorie compl\`ete de la construction et de la man\oe uvre des vaisseaux,  Saint Petersburg, Imprimerie de l'Acad\'emie Imp\'eriale ses Sciences, 1773. New augmented and corrected edition, ed. Claude-Antoine Jombert, fils a\^\i n\'e, Librairie du Roi pour le G\'enie et l'Artillerie, Paris, 1776.
Opera Omnia, Series 2, Volume 21, p. 80--222.

\bibitem{E520} L. Euler, Essai d'une th\'eorie de la r\'esistance qu'\'eprouve la proue d'un vaisseau dans son mouvement, M\'emoires de l'Acad\'emie des sciences de Paris, 1778, 1781, p. 597--602.
Opera Omnia, Series 2, Volume 21, p. 223--229 

\bibitem{Euler-Mditationes} L. Euler, Meditationes super problemate nautico, quod Illustrissima regia Parisienis Academia Scientiarum proposuit.  Opera omnia, Series 2, volume 20, p. 1--35.


\bibitem{Euler-Navalis} L. Euler, Scientia navalis seu tractatus de construendis ac dirigendis navibus, 2 volumes, Saint Petersburg, 1749. Opera Omnia, Series 2, vol. 18 and 19, Zurich and Basel, 1967 and 1972.

\bibitem{Euler-rep-1777} L. Euler, De repraesentatione superficiei sphaericae super plano, Acta Academiae Scientarum Imperialis Petropolitinae 1777, 1778, p. 107--132. 
Opera Omnia, Series 1, Volume 28, p. 248--275.


 \bibitem{Euler-pro-1777} L. Euler, De proiectione geographica superficiei sphaericae, Acta Academiae Scientarum Imperialis Petropolitinae 1777, 1778, p. 133-142.
Opera Omnia, Series 1, Volume 28, p. 276--287.
 
 
 \bibitem{Euler-pro-Desli-1777} L. Euler, De proiectione geographica Deslisliana in mappa generali imperii russici usitata, Acta Academiae Scientarum Imperialis Petropolitinae 1777, 1778, p. 143--153.
Opera Omnia, Series 1, Volume 28.



\bibitem{Euler-C} Leonhardi Euleri Commentationes arithmeticae collectae, 2 vols. Saint Petersburg, 1849.
          
          
\bibitem{Erm}     N. S. Ermolaeva, Mathematical cartography and D. A.  Grav\'e's method for solving the Dirichlet problem (Russian), Istor.-Mat. Issled. No. 32--33, (1990), 95--120.
  
  \bibitem{Gauss-Copenhagen} C. F. Gauss,  Allgemeine Aufl\"osung der Aufgabe: die Theile einer gegebnen Fl\"ache auf einer andern gegebnen Fl\"ache so abzubilden, da\ss \  die Abbildung dem Abgebildeten in den kleinsten Theilen \"ahnlich wird, Astronom. Abh. (Schumacher ed.), 3 (1825) 1--30. Also in Gauss's \emph{Werke}, vol. IV, 189-216. 




\bibitem{Gauss-English} K. F. Gauss,  \emph{Disquisitiones generales circa superficies curvas} (General investigations on curved surfaces), Gottingae, Typis Dieterichianis,1828.  Translated with notes and a bibliography by J. C. Morehead and A. M. Hilterbeitel, Princeton, 1902.



       \bibitem{Ghys}  E. Ghys, Sur la coupe des vêtements. Variation autour d'un thème de Tchebychev. Enseign. Math. 57 (2011), 165-208.
       
       \bibitem{Mes} G. Mescheryakou  The theoretical background of mathematical cartography (Russian), Moscow, 1968.
       
       \bibitem{Milnor} J. Milnor, A problem in cartography,
The American Mathematical Monthly
Vol. 76, No. 10 (Dec. 1969), 1101--1112.

\bibitem{Hilbert-Problems}   D. Hilbert,  Mathematische Probleme, \emph{G\"ottinger Nachrichten}, 1900, p. 253--297, reprinted in Archiv der Mathematik und Physik, 3d. ser., vol. 1 (1901) p. 44--63 and 213--237. English version, ``Mathematical problems", translated by M. Winston Newson, Bulletin of the AMS, vol. 8, 1902, p. 437-- 445 and 478--479. The English translation was also reprinted in  ``Mathematical developments arising from Hilbert problems", Proceedings of Symposia in Pure Math., Vol. XXVII, Part 1, F. Browder (Ed.), AMS, Providence, Rhode Island, 1974. Reprinted also in
 the  Bull. Amer. Math. Soc. (N.S.) 37 (2000), no. 4, 407-436.


\bibitem{H1900}   D. Hilbert,  \"Uber Flächen von konstanter Gau\ss scher Kr\"ummung. Transactions of the American mathematical Society 2,  (1901) 87--99. Also Appendix V to 
Hilbert's \emph{Foundations of geometry}, 2nd ed., Transl. L. Unger, Open court, La Salle, Illinois, 1971, 191--199.

\bibitem{Kempe} A. B. Kempe, On a General Method of describing Plane Curves of the nth degree by linkwork.
Proc. London Math. Soc. 7 (1876) 213--216.


 
\bibitem{Lagrange1779} J.-L. de Lagrange, 
Sur la construction des cartes g\'eographiques, Nouveaux m\'emoires de l'Acad\'emie royale des sciences et belles-lettres de Berlin, ann\'ee 1779, Premier m\'emoire,  \emph{\OE uvres compl\`etes}, tome 4,  637--664. Second m\'emoire,  \emph{\OE uvres compl\`etes}, tome 4, 664-692.  

 

\bibitem{Lambert-Bey}  J. H. Lambert, Beitr\"age zum Gebrauche der Mathematik und deren Anwendung, 3 parts,  im Verlage des Buchladens der Realschule,  Berlin,  1765-1772.


\bibitem{Lucas0} E. Lucas, Application de l'arithm\'etique \`a la construction de l'armure des satins r\'eguliers, Paris, G. Retaux, Nov. 1867. 16 pp.

\bibitem{Lucas1} E. Lucas, Lois g\'eom\'etriques des tissages, Congrès de l'Association Fran\c caise pour l'Avancement des Sciences  5 (1876), p. 114.

\bibitem{Lucas2} E. Lucas, Les principes fondamentaux de la g\'eom\'etrie des tissus, Congrès de l'Association Fran\c caise pour l'Avancement des Sciences (2), 40 (1911), p. 72-88, m\'emoire extrait de l'\emph{Ingegnere Civile}, Turin, 1880, traduit de l'Italien et condens\'e par A. Aubry et A. G\'erardin.


    \bibitem{Papa-Inde2} A. Papadopoulos, On the works of Euler and his followers on spherical geometry,  \emph{Ga\d{n}ita Bh\=ar\=at\=\i } (Indian Mathematics), the Bulletin of the Indian Society for History of Mathematics, 36 No. 1-2,  (2014) 237-292.

\bibitem{Papa-qc} A. Papadopoulos, Quasiconformal mappings,  
 from Ptolemy's geography to the work of Teichm\"uller, to appear in: \emph{Uniformization, Riemann-Hilbert Correspondence, Calabi-Yau Manifolds, and Picard-Fuchs Equations} (ed. L. Ji and S.-T. Yau),  International Press and  Higher Education Press. To appear in 2017.

\bibitem{bio} C. A. Poss\'e, Excerpts of a biography of Chebyshev, contained in his \emph{Collected Works}, Edited by A. Markoff and N. Sonin, Vol. II,  p. I-VI.


\bibitem{Rosenfeld} B. A. Rosenfeld and N. E. Maryukova, Surfaces of constant curvature and geometric interpretations of the Klein-Gordon, sine-Gordon equations, Publications de l'Institut Math\'ematique (Beograd) Nouvelle S\'erie, tome 61 (75) (1997), 119--132.

 \bibitem{S0} F. T. von Schubert, De proiectione Sphaeroidis ellipticae geographica, Dissertatio Prima, \emph{Nova Acta Academiae Scientiarum Imperialis petropolitana}, vol. V,  1789, p. 130-146.
 
   \bibitem{S00} F. T. von Schubert, De proiectione Sphaeroidis ellipticae geographica, Dissertatio Secunda. Proiectio stereographica horizontalis feu obliqua, \emph{Nova Acta Academiae Scientiarum Imperialis petropolitana}, vol. VI,   1790, p. 123-142.
   
     \bibitem{S01} F. T. von Schubert, De proiectione Sphaeroidis ellipticae geographica, Dissertatio Tertia. \emph{Nova Acta Academiae Scientiarum Imperialis petropolitana}, vol. VII,   1791, p. 149-161.
  



\bibitem{Sossinsky} A. Sossinsky,   Configuration spaces of planar linkages, In: Handbook of Teichm\"uller theory (ed. A. Papadopoulos) Vol. VI, p. 335--373, European Mathematical Society, Zurich, 2016.



 \bibitem{Thurston1985} W. P.
Thurston, \emph{Minimal Stretch maps between hyperbolic
surfaces}, Preprint (1985) available at
http://arxiv.org/abs/math.GT/9801039



\bibitem{TW} W. P. Thurston and J. Weeks, The Mathematics of three dimensional manifolds, \emph{Scientific American}, 251 (1984) 94-106. 

\bibitem{TWF} W. P. Thurston and J. Weeks,  Les vari\'et\'es \`a  trois dimensions, \emph{Pour la science}, 84, October 1984. 




\bibitem{TP} W. P. Thurston, \emph{Three-dimensional geometry and topology} Princeton University Press, Princeton, 1997.

    \bibitem{Ti} V. M. Tikhomirov, Approximation theory in the twentieth century, In: 
{Mathematical Events of the Twentieth Century}, p. 409-436 Springer Verlag, Berlin, 2006.


\bibitem{Wassilief} A. Vassilief, P. Tch\'ebychef et son \oe uvre scientifique, Bollettino di bibliografia e storia delle scienze matematiche, I Turin, 1898. 

\bibitem{You} A. P. Youschkevitch, Chebyshev, Pafnuty Lvovich, In: Complete Dictionary of Scientific Biography, Charles Scribner's Sons, 2008.

\end{thebibliography}
\end{document}